\documentclass[12pt]{amsart}
\usepackage{amsmath,amsthm, amscd, amssymb, amsfonts}

\newcommand{\C}{{\mathcal C}}

\newcommand{\Ind}{\mbox{\rm Ind\,}}
\newcommand\op{\operatorname{op}}
\newcommand\ev{\operatorname{ev}}
\newcommand\coev{\operatorname{coev}}

\newcommand\vect{\operatorname{Vec}}

\newcommand\Rep{\operatorname{Rep}}

\newcommand\Hom{\operatorname{Hom}}
\newcommand\Opext{\operatorname{Opext}}

\newcommand{\fde}{{\triangleright}}
\newcommand{\fiz}{{\triangleleft}}

\numberwithin{equation}{section}\theoremstyle{plain}

\newtheorem{theorem}{Theorem}[section]
\newtheorem{lema}[theorem]{Lemma}
\newtheorem{cor}[theorem]{Corollary}

\newtheorem{prop}[theorem]{Proposition}
\newtheorem{claim}{Claim}[section]

\theoremstyle{definition}
\newtheorem{definition}[theorem]{Definition}

\newtheorem{question}[equation]{Question}

\theoremstyle{remark}
\newtheorem{obs}[theorem]{Remark}

\newcommand\id{\operatorname{id}}

\def\pf{\begin{proof}}
\def\epf{\end{proof}}

\theoremstyle{remark}

\begin{document}

\renewcommand{\baselinestretch}{1.2}
\renewcommand{\thefootnote}{}
\thispagestyle{empty}
\title[Frobenius-Schur indicators]{Frobenius-Schur indicators for a class of fusion categories}
\author{Sonia Natale}
\address{Facultad de Matem\'atica, Astronom\'\i a y F\'\i sica
\newline \indent
Universidad Nacional de C\'ordoba
\newline
\indent CIEM -- CONICET
\newline
\indent (5000) Ciudad Universitaria
\newline
\indent
C\'ordoba, Argentina}
\email{natale@mate.uncor.edu \newline
\indent \emph{URL:}\/ http://www.mate.uncor.edu/natale}
\thanks{This work was partially supported by CONICET,
CONICOR, Fundaci\' on Antorchas  and Secyt (UNC)}
\date{Revised version of October 5, 2004.}

\begin{abstract} We give an explicit description, up to gauge equivalence, of group-theoretical quasi-Hopf algebras.
We use this description to compute the Frobenius-Schur indicators for  group-theoretical fusion categories.
\end{abstract}

\maketitle

\section{Introduction}

There are strong analogies between the theory of finite groups and
the theory of semisimple Hopf algebras; some of them, however,
still remain conjectural. In particular, the problem of
classifying semisimple Hopf algebras, say over the field of
complex numbers, seems to be a considerably difficult one, even in
low dimensions. Perhaps the most important feature of these
objects, which relates them to other branches of mathematics and
physics, is that their category of representations is a special
case of a so called \emph{fusion category}. This fact leads to the
consideration of the classification problem, not only modulo Hopf
algebra isomorphisms, but modulo \emph{gauge equivalences}:
roughly,  two finite dimensional (quasi)-Hopf algebras $H$ and
$H'$ give rise to the same fusion category of representations if
and only if they are \emph{gauge equivalent}, in the sense that $H
= H'$ as algebras, and  the comultiplication of $H'$  is obtained
by 'twisting' that of $H$ by means of $\Delta_{H'}(h) = F
\Delta(h)F^{-1}$, for some \emph{gauge transformation} $F \in (H
\otimes H)^{\times}$.

\bigbreak An  important class of examples of semisimple quasi-Hopf
algebras was introduced by Ostrik \cite{ostrik} and studied later
by Etingof, Nikshych and Ostrik \cite{ENO}: these are called
\emph{group theoretical} and, by definition, they are exactly
those for which the category of representations is a group
theoretical  category $\C(G, \omega, F, \alpha)$, where $G$ is a
finite group, $F \subseteq G$ is a subgroup,   $\omega: G \times G
\times G \to k^{\times}$ is a normalized 3-cocycle and  $\alpha: F
\times F \to k^{\times}$ is a normalized 2-cochain, such that
$\omega\vert_F = d\alpha$. More precisely, $\C(G, \omega, F,
\alpha)$ is the category of $k_{\alpha}F$-bimodules in the tensor
category $\vect^G_{\omega}$ of finite dimensional $G$-graded
vector spaces, with associativity constraint given by $\omega$.
All concrete known examples of semisimple Hopf algebras (and this
does not extend to the \emph{quasi}-context), turn out to be group
theoretical;  this can be seen as a consequence of our results in
\cite{gp-ttic}. The following question  has been posed in
\cite{ENO}. An answer to this question (even an affirmative one)
would be of great significance in the classification program.

\begin{question} Does there exist a semisimple Hopf algebra which is not group theoretical?
\end{question}

The categorical nature of this question leads naturally to the problem of finding and computing  \emph{gauge  invariants} of group theoretical quasi-Hopf algebras, that is, invariants which depend on the gauge equivalence class of the object rather than on the isomorphism class itself.

\bigbreak Recently, Mason and Ng have constructed a gauge
invariant,  the \emph{Frobenius-Schur} indicators, for semisimple
quasi-Hopf algebras \cite{m-ng}. They have proved {\it loc. cit.}
a generalization of the Frobenius-Schur Theorem for finite groups,
c.f. \cite{serre}. Their construction extends results of Linchenko
and Montgomery for semisimple Hopf algebras \cite{li-mo}; it also
extends results of Bantay on the Frobenius-Schur indicators for
the Dijkgraaf-Pasquier-Roche quasi-Hopf algebra $D^{\omega}G$,
after the definition in \cite{bantay} of the indicators attached
to conformal field theories.

Essentially, Frobenius-Schur indicators were  defined in a
categorical fashion for any semisimple rigid tensor category which
is \emph{pivotal}, \textit{i.e.}, which admits a natural tensor
isomorphism between the identity and the second left duality
functors, in the work of Fuchs, Ganchev, Szlach\' anyi and
Vescerny\' es \cite{fgsv}. It is shown in \cite{ENO} that
representation categories of semisimple quasi-Hopf algebras are in
fact pivotal.

Other gauge invariants can be attached to a semisimple quasi-Hopf
algebra. One of the most studied is the $K_0$-ring of its
representation category. This invariant does not distinguish the
group algebras of the two nonabelian groups of order 8: the
dihedral group $D_4$ and the quaternionic group $Q_2$. However, by
a result of Tambara and Yamagami \cite{TY}, these two groups are
not gauge equivalent. In their paper, Mason and Ng have noted that
the Frobenius-Schur indicators do distinguish the dihedral and
quaternionic groups. As pointed out to us by the referee, in some
cases it may happen that the Frobenius-Schur indicators contain
less information than the $K_0$-ring, e.g., in the case of dual
group algebras. So, in some sense, these two invariants are of
very different nature.

\bigbreak
In this paper we give an explicit description of group-theoretical Hopf algebras and use it to compute their irreducible characters and Frobenius-Schur indicators. Our main original contributions in the  description of the quasi-structure are on the one hand the proof of the existence of a certain normalization of the 3-cocycle $\omega$ (Proposition \ref{ln-om}), and on the other hand,  the construction of a quasi-antipode (Theorem \ref{q-antip}) which is, of course,  esential in the computation of the Frobenius-Schur invariants.
We obtain the following formula for the indicator of the irreducible character $\chi$:
\begin{equation}\chi(\nu_{A^{\op}})  =
\vert F \vert^{-1}   \sum_{(q \fiz x).q = e} \omega(xq, xq, xq)  \quad \chi(\delta_{q} (xq)^2). \end{equation}
which involves a certain normalization of the 3-cocycle $\omega$. See Corollary \ref{fs-chi}. Here, $q$ runs over an appropriate choice of representatives of $G$ modulo $F$.

One instance of these examples comes from an exact factorization $G = F Q$ of the group $G$ into its subgroups $F$ and $Q$.  In this case, there is a group $\Opext(kF, k^Q)$, which classifies the \emph{abelian} Hopf algebra extensions of $k^Q$ by $kF$: as a Hopf algebra, the extension corresponding to an element $[\sigma, \tau] \in \Opext(kF, k^Q)$ is a  bicrossed product $k^Q \#_{\sigma}^{\tau}kF$, where $\sigma: F \times F \to (k^Q)^{\times}$ and $\tau: Q \times Q \to (k^F)^{\times}$ are a pair of compatible cocycles; see \cite[Theorem 1.2]{gp-ttic}. In this case, the 3-cocycle $\omega$ in our formula for the Frobenius-Schur indicators is the one associated to $[\sigma, \tau]$ in the Kac exact sequence.  This gives an alternative compact expression for the formula found by Kashina, Mason and Montgomery in \cite{KMM}.

We would like to point out that the description for the quasi-Hopf algebra structure for group-theoretical quasi-Hopf algebras generalizes the construction of the twisted quantum doubles $D^{\omega}G$ by Dijkgraaf, Pasquier and Roche. This agrees with the characterization  given in our paper \cite{gp-ttic} in terms of quantum (or Drinfeld) doubles; so in some sense  these quasi-Hopf algebras are all of \emph{DPR-type}.

\bigbreak
The paper is organized as follows. In Sections \ref{fs-mng} and \ref{gt-fc} we recall the definition  of the indicators constructed by Mason and Ng \cite{m-ng} and the definition and main properties  of group-theoretical categories as given by Etingof, Nikshych  and Ostrik \cite{ostrik}, \cite{ENO}. In Section \ref{qha} we give a description, up to gauge equivalence, of the structure of group-theoretical quasi-Hopf algebras, and finally in Section \ref{fs} we give an explicit formula for the Frobenius-Schur indicators of group theoretical categories.
We consider some examples in Section \ref{ejemplos}.
Throughout this paper we work over an  algebraically closed field $k$ of
characteristic zero.

\section{Frobenius-Schur indicators}\label{fs-mng}

 We shall recall the definition  of the indicators constructed by Mason and Ng \cite{m-ng}.

\subsection{} Let $(H, \Delta, \epsilon, \phi, \mathcal S, \alpha, \beta)$ be a finite dimensional semisimple quasi-Hopf algebra \cite{drinfeld} (later on indicated $(H, \Phi)$ for short),
 that is,  $H$ is an associative unital algebra over $k$ which is semsimple and finite-dimensional; $\epsilon: H \to k$ and $\Delta: H \to H \otimes H$ are algebra maps; $\Phi \in H^{\otimes 3}$ is an invertible element such that
\begin{equation}\label{f1}(\id \otimes \id \otimes \Delta)(\Phi) (\Delta \otimes \id \otimes \id)(\Phi) = (1 \otimes \Phi) (\id \otimes \Delta \otimes \id)(\Phi) (\Phi \otimes 1), \end{equation}
\begin{equation}\label{f2}(\id \otimes \epsilon \otimes \id)(\Phi)= 1 \otimes 1, \end{equation}
\begin{equation}\label{f3}(\epsilon \otimes \id)\Delta (h)  = h = (\id \otimes \epsilon)\Delta (h), \end{equation}
\begin{equation}\label{f4}\Phi (\Delta \otimes \id)\Delta (h) \Phi^{-1}  = (\id \otimes \Delta)\Delta (h), \end{equation} for all $h\in H$.
The map $\mathcal S: H \to H^{\op}$ is an algebra anti-automorphism of $H$; $\alpha, \beta \in H$ are such that
\begin{equation}\label{ant-ab}\mathcal S(h_1) \alpha h_2 = \epsilon(h)\alpha, \qquad h_1 \beta \mathcal S(h_2) = \epsilon (h) \beta, \qquad \forall h\in H; \end{equation}
\begin{equation}\label{largas}\Phi^{(1)} \beta \mathcal S(\Phi^{(2)}) \alpha \Phi^{(3)} = 1 = \mathcal S(\Phi^{(-1)}) \alpha \Phi^{(-2)} \beta \mathcal S(\Phi^{(-3)}),
 \end{equation}
where we are using the abbreviated notation $\Phi = \Phi^{(1)} \otimes \Phi^{(2)} \otimes \Phi^{(3)}$ and $\Phi^{-1} = \Phi^{(-1)} \otimes \Phi^{(-2)} \otimes \Phi^{(-3)}$.

\bigbreak
The category $\Rep H = : \Rep (H, \Phi)$ is a \emph{fusion category}, in the terminology of \cite{ENO}. The associativity constraint is given by the natural action of $\Phi$; the left  dual of an object $V$ of $\Rep H$ is the vector space $V^* = \Hom (V, k)$ with the $H$-action $\langle h.f, v \rangle = \langle f, \mathcal S(h)v \rangle$; and the evaluation and coevaluation maps are given, respectively, by
\begin{equation}\ev: V^* \otimes V \to k, \qquad  \ev(f \otimes v) = \langle f, \alpha . v \rangle, \end{equation}
\begin{equation} \coev: k \to  V \otimes V^*, \qquad  1 \mapsto \sum_i \beta . v_i \otimes  v^i,\end{equation} for all $f \in V^*$, $v\in V$, where $(v_i)$ and $(v^i)$ are dual basis of $V$.

\bigbreak
Two quasi-Hopf algebras $H_1$ and $H_2$ are called \emph{gauge equivalent} if there exists a \emph{gauge transformation},  {\it i.e.}, an invertible  normalized element  $F \in H_1 \otimes H_1$ such that $(H_1)_F$ and $H_2$ are isomorphic as quasi-bialgebras.

Here, $(H_1)_F$ is the quasi-Hopf algebra $(H_1, \Delta_F, \epsilon, \Phi_F, \mathcal S_F, \alpha_F, \beta_F)$, where
\begin{equation*}\Delta_F(h) = F \Delta(h)F^{-1},  \qquad h \in H, \end{equation*}
\begin{equation*}\Phi_F = (1 \otimes F) (\id \otimes \Delta)(F) \Phi (\Delta \otimes \id)(F^{-1}) (F^{-1} \otimes 1),  \end{equation*}
\begin{equation*}\alpha_F = \mathcal S(F^{(-1)}) \alpha F^{(-2)},  \qquad \beta_F = F^{(1)} \beta \mathcal S(F^{(2)}); \end{equation*}
where $F = F^{(1)}  \otimes F^{(2)}$, $F^{-1} = F^{(-1)} \otimes F^{(-2)}$.

\bigbreak
The finite dimensional quasi-Hopf algebras $H_1$ and $H_2$ are gauge equivalent if and only if $\Rep H_1$ is equivalent to $\Rep H_2$ as $k$-linear tensor categories. See  \cite{et-gel}.

\begin{obs} It is shown in \cite{ENO} that the fusion categories of the form $\Rep (H, \Phi)$ are exactly those for which the Frobenius-Perron dimensions of simple objects are integers. \end{obs}

\subsection{} Let $(H, \Phi)$ be a finite dimensional quasi-Hopf algebra. A \emph{normalized two sided integral} of $H$ is an element $\Lambda \in H$ such that $$h\Lambda = \epsilon(h)\Lambda = \Lambda h, \quad \forall h \in H; \qquad \epsilon(\Lambda) = 1.$$
Suppose that $(H, \Phi)$ is semisimple. Then $H$ contains a unique normalized two sided integral \cite{h-nill}.

The following definition is due to Mason and  Ng \cite{m-ng}. It generalizes a previous definition for semisimple Hopf algebras given by  Linchenko and Montgomery \cite{li-mo}.

\begin{definition} Let $(H, \Phi)$ be a finite dimensional semisimple quasi-Hopf algebra and let $\Lambda \in H$ be a normalized two sided integral. Let also $\chi \in H^*$ be an irreducible character of $H$. The \emph{Frobenius-Schur indicator} of $\chi$ is the element $\chi(\nu_H)$, where
$\nu_H$ is the canonical central element of $H$ given by
\begin{equation}\nu_H = m (q_L \Delta(\Lambda) p_L); \end{equation}
here, $m: H \otimes H \to H$ is the multiplication map, and $q_L, p_L \in H^{\otimes 2}$ are defined by
\begin{equation*}q_L : = \mathcal S(\Phi^{(-1)}) \alpha \Phi^{(-2)} \otimes \Phi^{(-3)}, \qquad p_L : = \Phi^{(2)} \mathcal S^{-1}(\Phi^{(1)} \beta) \otimes \Phi^{(3)}. \end{equation*}
\end{definition}

The family of Frobenius-Schur indicators $\{ \chi(\nu_H) \}_{\chi}$ is an invariant of the $k$-linear tensor category $\Rep (H, \Phi)$. This means that it is invariant under gauge transformations of $(H, \Phi)$.

Also, if $\alpha$ and $\beta$ are invertible elements of $H$, then the canonical central element $\nu_H$ can be computed as follows \cite[Corollary 3.5]{m-ng}:
\begin{equation}\nu_H = (\Lambda_1\Lambda_2) (\beta \alpha)^{-1} = (\beta \alpha)^{-1} (\Lambda_1\Lambda_2).
\end{equation}

\bigbreak
In analogy with finite group situation, the Frobenius-Schur indicator of the irreducible character $\chi = \chi_V$ satisfies the following:

(i) $\chi(\nu_H) = 0,  1$  or $-1$,  and $\chi(\nu_H) \neq 0$  if and only if $\chi = \chi^*$

(ii) $\chi(\nu_H) = 1$ (respectively $-1$) if and only if $V$  admits a non-degenerate bilinear form $\langle \, , \, \rangle: V \otimes V \to k$, with adjoint $\mathcal S$,
such that $\langle x, y \rangle = \langle y, g^{-1}x \rangle$ (respectively, $\langle x, y \rangle = - \langle y, g^{-1}x \rangle$),  where $g \in H$ is the so called \emph{trace element} of $H$.

\section{Group theoretical fusion categories}\label{gt-fc}

Group theoretical categories were introduced in \cite[Section 3]{ostrik} and also studied in \cite{ENO}. In this section we recall their definition and basic properties.

\subsection{}\label{gt-cat} Let $G$ be a finite group, and let $F \subseteq G$ be a subgroup. The identity element of $G$ will be denoted by $e$. Let also be given the following data:

\begin{itemize}\item a normalized 3-cocycle $\omega: G \times G \times G \to k^{\times}$, that is,
\begin{align}\label{cc-omega} \omega(ab, c, d) \omega(a, b, cd) & = \omega(a, b, c) \omega(a, bc, d) \omega(b, c, d), \\
\omega(e, a, b) & = \omega(a, e, b) = \omega(a, b, e) = 1,\end{align}
for all $a, b, c, d \in G$;
\item a normalized 2-cochain $\alpha: F \times F \to k^{\times}$;
\end{itemize}
subject to the condition
\begin{equation}\label{rest-to-F}\omega\vert_{F\times F \times F} = d\alpha. \end{equation}

Consider the category
$\vect^G_{\omega}$ of finite dimensional $G$-graded vector spaces, with
associativity constraint given by $\omega$: explicitly, for any
three objects $U, U'$ and $U''$ of $\vect^G_{\omega}$, we
have $a_{U, U', U''}: (U \otimes U') \otimes U'' \to U \otimes (U'
\otimes U'')$, given by
\begin{equation}
a_{U, U', U''} ((u \otimes u') \otimes u'') = \omega(||u||,
||u'||, ||u''||) \, u \otimes (u' \otimes u''),
\end{equation}
on homogeneous elements $u \in U$, $u' \in U'$, $u'' \in U''$,
where we use the symbol $|| \, ||$ to denote the corresponding
degree of homogeneity. In other words, $\vect^G_{\omega}$ is
the category of representations of the quasi-Hopf algebra
$k^{G}$, with associator $\omega \in (k^{G})^{\otimes 3}$.

By \eqref{rest-to-F}, the twisted group algebra $k_{\alpha}F$ is an (associative unital) algebra in $\vect^G_{\omega}$, and one may naturally attach to it a monoidal category.
Precisely, the category $\C(G, \omega, F, \alpha)$ is by definition the $k$-linear monoidal category  of $k_{\alpha}F$-bimodules in $\vect^G_{\omega}$: tensor product is $\otimes_{k_{\alpha}F}$ and the unit object is $k_{\alpha}F$. This is a fusion category over $k$ with the property that the Frobenius-Perron dimensions of its objects are integers \cite[8.8]{ENO}.

\bigbreak
The categories of the form $\C(G, \omega, F, \alpha)$ are called \emph{group theoretical} \cite[Definition 8.46]{ENO}.
By extension, a (quasi)-Hopf algebra $A$ is called group theoretical if the category $\Rep A$ of its finite dimensional representations is group theoretical.

\subsection{} Let $\eta: G \times G \to k^{\times}$ and $\chi: F \to k^{\times}$ be normalized cochains, and let $\widetilde\omega: G \times G \times G \to k^{\times}$, $\widetilde \alpha: F \times F \to k^{\times}$ be given by
\begin{equation}\widetilde \omega = \omega (d\eta), \qquad \widetilde \alpha = \alpha (\eta\vert_{F\times F})  (d\chi).\end{equation}
 Then the categories $\C(G, \omega, F, \alpha)$ and $\C(G, \widetilde \omega, F, \widetilde \alpha)$ are equivalent \cite[Remark 8.39]{ENO}.

\begin{obs}\label{alfa=1} Let $G$, $F$, $\omega$ and $\alpha$ be as above. Let $Q$ be a set of representatives of the left cosets of $F$ in $G$ such that $e \in Q$; so that every element $g \in G$ writes uniquely in the form $g =  xp$, with $p \in Q$, $x \in F$.
Consider the 2-cochain $\eta: G \times G \to k^{\times}$ defined in the form
\begin{equation} \eta(xp, yq): = \alpha^{-1}(x, y), \qquad p, q \in Q, \quad x, y \in F.\end{equation}
Then, taking $\chi =  1$, we obtain $\widetilde \alpha = 1$. Therefore the categories $\C(G, \omega, F, \alpha)$ and $\C(G, \widetilde \omega, F, 1)$ are equivalent, where $\widetilde \omega =  \omega (d \eta)$. That is, up to monoidal equivalence, \emph{we may always assume that $\alpha = 1$.}

\medbreak
Note also that the categories $\C(G, \omega, F, 1)$ and $\C(G, \omega (d \eta), F, 1)$ are tensor equivalent for every normalized 2-cochain $\eta: G \times G \to k^{\times}$ such that $\eta\vert_{F \times F}$ is a \emph{coboundary}.
\end{obs}

\subsection{}The fiber functors $\C(G, \omega, F, \alpha) \to \vect$, in the case they exist, are classified by conjugacy classes of subgroups $\Gamma$ of $G$, endowed with a
2-cocycle $\beta \in Z^2(\Gamma,
k^{\times})$,  such that the class of $\omega\vert_{\Gamma}$ is trivial; $G = F \Gamma$ and the class of the cocycle $\alpha\vert_{F \cap \Gamma} \beta^{-1}\vert_{F \cap
\Gamma}$ is non-degenerate \cite[Corollary 3.1]{ostrik}.

\begin{obs} The category $\C = \C(G, \omega, F, \alpha)$ has the property that the Frobenius-Perron dimensions of its objects are integers. A Tannaka-Krein reconstruction argument shows,  that $\C$ is equivalent to the category of representations of a semisimple quasi-Hopf algebra over $k$ \cite[Theorem 8.33]{ENO}.  \end{obs}

\bigbreak
It follows from \cite[8.8]{ENO} that duals, opposites, quotient categories, full subcategories, and tensor products of group theoretical categories are also group theoretical. Also, by \cite[Remark
8.47]{ENO},  the Drinfeld center $\mathcal Z(\C)$ is group theoretical if and only if so is $\C$.

However, in Remark 8.48 of the paper \cite{ENO}, the authors note that there exist semisimple quasi-Hopf algebras such that their category of representations are \emph{not} group theroretical: an explicit example is quoted in {\it loc. cit.} which comes from the construction of Tambara and Yamagami \cite{TY}. The answer to the corresponding question for semisimple Hopf algebras is still not known.

\section{Group theoretical quasi-Hopf algebras}\label{qha}

The aim of this section is to give an explicit description, up to gauge equivalence, of the structure of group-theoretical quasi-Hopf algebras.
This will enable us to explicitly compute the Frobenius-Schur indicators of group theoretical categories in the next section.
The description is based on a result of Schauenburg \cite[3.4]{sb}, which reconstructs a quasi-bialgebra structure from certain monoidal categories of bimodules in a more general context.

Our main new result concerning this description is the explicit construction of the quasi-antipode in the group theoretical case, which is relevant for our purposes; see Theorem \ref{q-antip}.

An instance of this quasi-Hopf algebra construction, for the case where $\omega = 1$ and $\alpha = 1$, was studied by Y. Zhu in \cite{zhu}.  This case was  also studied in \cite{beggs, beggs-mod}, from the point of view of  the tensor categories of representations. Throughout this paper we shall adopt the notation in \cite{beggs}.

\bigbreak
In what follows we shall fix a finite group $G$ and a subgroup $F \subseteq G$. Following \cite{zhu}, we shall also fix a set of simultaneous representatives of the left and right cosets of $F$ in $G$, $Q \subseteq G$; this is possible since $G$ is finite. Thus every element $g \in G$ has  unique factorizations $g = xq = py$, where $x, y \in F$, $q, p  \in Q$. We assume that $e \in Q$.

\subsection{}
The uniqueness of the factorization $G = FQ$ implies that there are well defined maps $$\fde: Q \times F \to F, \qquad  \fiz: Q \times F \to Q,$$
$$.: Q \times Q \to Q, \qquad \theta: Q \times Q \to F,$$ determined by the conditions
\begin{align}qx & = (q \fde x) (q \fiz x), \qquad q \in Q, \, x \in F;\\
pq & = \theta(p, q) p.q, \qquad p, q \in Q.\end{align}

The main relations between these maps are stated in the following lemma.

\begin{lema}\label{prop-beggs} (\cite[Proposition 2.4]{beggs}.) The following identities hold, for all $p, q, r \in Q$, $x, y \in G$:

\medbreak
\noindent (i) $p \fiz xy = (p \fiz x) \fiz y$, $p \fiz e = p$;

\medbreak
\noindent (ii) $(p.q) \fiz x = (p \fiz (q \fde x)) . (q \fiz x)$;

\medbreak
\noindent (iii) $p \fde (q \fde x) = \theta(p, q) \left( (p . q) \fde x \right) \theta((p \fiz (q \fde x), q \fiz x)^{-1}$, $e \fde x = x$;

\medbreak
\noindent (iv) $p \fde xy = (p \fde x)((p \fiz x) \fde y)$;

\medbreak
\noindent (v) $\theta(p, q) \theta(p.q, r) = (p \fde \theta(q, r)) \theta(p \fiz \theta(q, r), q.r)$;

\medbreak
\noindent (vi) $(p \fiz \theta(q, r)) . (q.r) = (p.q).r$.

\medbreak
\noindent (vii) $\theta(p, e) = \theta(e, p) = e$.\qed \end{lema}

\subsection{} Let $\omega: G \times G \times G \to k^{\times}$ be a normalized 3-cocycle such that $\omega\vert_{F\times F \times F}$ is trivial.
In what follows we shall fix the group theoretical category $\C = \C(G, \omega, F, 1)$.

Thus the cochain $\alpha: F \times F \to k^{\times}$ as in Subsection \ref{gt-cat} will be trivial. This is, up to monoidal equivalence, no loss of generality thanks to Remark \ref{alfa=1}.

\begin{prop}\label{ln-om} There exists a normalized 2-cochain $\eta: G \times G \to k^{\times}$ such that $\eta\vert_{F \times F} = 1$ and $\omega(d\eta)\vert_{F \times G \times G} = 1 = \omega(d\eta)\vert_{F \times F \times Q}$. \end{prop}

\pf Recall that the coboundary $d\eta: G \times G \times G \to k^{\times}$ is given by
$$(d\eta) (a, b, c) = \eta(ab, c) \,  \eta(a, b) \, \eta(b, c)^{-1} \,  \eta(a, bc)^{-1},$$
for all $a, b, c \in G$.

The proof will be done in three steps. Let first $\eta_1: G \times G \to k^{\times}$ be the normalized cochain given by
$$\eta_1(xp, yq): = \omega(x, y, q), \qquad x, y \in F, \quad p, q \in Q.$$
Then we have $\eta_1\vert_{G \times F} = 1$ and for all $x, y, z \in F$, $q \in Q$, we have
\begin{align*}(d\eta_1)(x, y, zq) & = \eta_1(xy, zq) \, \eta_1(x, y) \, \eta_1(y, zq)^{-1} \, \eta_1(x, yzq)^{-1} \\
& = \omega(xy, z, q) \, \omega(y, z, q)^{-1} \, \omega(x, yz, q)^{-1} \, \omega(x, y, z)^{-1} \\
& = \omega(x, y, zq)^{-1},
\end{align*}
the second equality because $\omega\vert_{F \times F \times F} = 1$. Thus $\omega(d\eta_1)\vert_{F \times F \times G} = 1$.

\bigbreak
Put now $\omega_0 = \omega(d\eta_1)$ and define $\eta_2: G \times G \to k^{\times}$ in the form
$$\eta_2(xp, yq): = \omega_0(x, p, yq) \, \omega_0(p, y, q)^{-1}, \qquad x, y \in F, \quad p, q \in Q.$$
Then $\eta_2\vert_{F \times G} = 1$ and we have
\begin{align*}(d\eta_2)(x, yp, zq) & = \eta_2(xyp, zq) \, \eta_2(x, yp)\, \eta_2(yp, zq)^{-1} \, \eta_2(x, ypzq)^{-1} \\
& = \omega_0(xy, p, zq) \, \omega_0(p, z, q) \, \omega_0(p, z, q)^{-1} \, \omega_0(y, p, zq)^{-1} \\
& = \omega_0(x, yp, zq) \, \omega_0(x, y, p) \, \omega_0(x, y, pzq)^{-1} \\
& = \omega_0(x, yp, zq),
\end{align*}
for all $x, y, z \in F$, $p, q \in Q$,
where in the third and fourth equalities we have used that $\omega_0\vert_{F \times F \times G} = 1$.
Hence $\omega_0(d\eta_2^{-1})\vert_{F \times G \times G} = 1$.

\bigbreak
Finally, let $\omega_1 = \omega_0(d\eta_2^{-1})$. The condition $\omega_1\vert_{F \times G \times G} = 1$ is equivalent to
$\omega_1(xt, g, h) = \omega_1(t, g, h)$,
for all $x \in F$, $t, g, h \in G$. Hence
$$\omega_1(zp, x, yq) = \omega_1(p, x, y) \, \omega_1(p, xy, q) \, \omega_1(px, y, q)^{-1},$$ for all $z \in F$.

Let $\eta_3: G \times G \to k^{\times}$ be defined by
$$\eta_3(xp, yq): = \omega_1(xp, y, q), \qquad x, y \in F, \quad p, q \in Q.$$
Then $\eta_3\vert_{F \times G} = \eta_3\vert_{G \times F}  = 1$, and
for all $x, y \in F$, $p, q \in Q$,
\begin{align*}(d\eta_3)(p, x, yq) & = \eta_3(px, yq) \, \eta_3(p, x)\, \eta_3(x, yq)^{-1} \, \eta_3(p, xyq)^{-1} \\
& = \eta_3(px, yq) \, \eta_3(p, xyq)^{-1} \\
& = \omega_1(px, y, q) \, \omega_1(p, xy, q)^{-1} \\
& =  \omega_1(p, x, y)\, \omega_1(p, x, yq)^{-1}.
\end{align*}
In particular, $(d\eta_3)(p, x, y) = 1$, and thus \begin{equation}\label{aux}\omega_1(d\eta_3)(p, x, yq) = \omega_1(d\eta_3)(p, x, y).\end{equation}

\begin{claim}We have $\omega_1(d\eta_3)\vert_{F \times G \times G} = 1$. \end{claim}

\pf Let $x, y, z \in F$, $p, q \in Q$. Using that $\omega_1\vert_{F \times G \times G} = 1$, we compute
\begin{align*}(d\eta_3)(x, yp, zq) & = \eta_3(xyp, zq) \, \eta_3(x, yp)\, \eta_3(yp, zq)^{-1} \, \eta_3(x, ypzq)^{-1} \\
& = \omega_1(xyp, z, q) \, \omega_1(yp, z, q)^{-1} \\
& =  \omega_1(p, z, q) \, \omega_1(p, z, q)^{-1} = 1.
\end{align*}
This proves the claim. \epf

In view of the claim, equation \eqref{aux} is equivalent to  $\omega_1(d\eta_3)\vert_{G \times F \times Q} = 1$.
This implies the proposition, since by construction $\omega_1(d\eta_3) = \omega (d\eta)$, for a suitable normalized 2-cochain such that $\eta\vert_{F \times F} = 1$. \epf

By Remark \ref{alfa=1}, the property $\eta\vert_{F \times F} = 1$ in Proposition \ref{ln-om} implies that $\C(G, \omega, F, 1)$ is  tensor equivalent to $\C(G, \omega (d\eta), F, 1)$. \emph{Then we may and shall assume in what follows that the 3-cocycle $\omega: G \times G \times G \to k^{\times}$ satisfies the normalization conditions}
\begin{equation}\label{norm-om}\omega\vert_{F\times G \times G} = 1,\end{equation}
\begin{equation}\label{norm-om-q}\omega\vert_{G\times F \times Q} = 1.\end{equation}

\bigbreak
These conditions are necessary in order to apply the results of \cite[3.4]{sb}; see Definition 3.3.2 in {\it loc. cit.}

\begin{lema}\label{eq-om} Let $g, h \in G$, $x, y \in F$, $p, q \in Q$. Then we have

(i) $\omega(xp, g, h) = \omega(p, g, h)$;

(ii) $\omega(g, y, xp) = \omega(g, y, x)$;

(iii) $\omega(g, x, pq) = \omega(g, x, \theta(p, q))$;

(iv) $\omega(pq, g, h) = \omega(p.q, g, h)$. \end{lema}

\pf Parts (i) and (ii) follow from the cocycle condition \eqref{cc-omega} and the normalization conditions \eqref{norm-om} and \eqref{norm-om-q}. Parts (iii) and (iv) are a consequence of parts (i) and (ii), respectively. \epf

\subsection{}\label{str-qha}
Let $A = k^Q \#_{\sigma} kF$ be the crossed product corresponding to the action $\rightharpoonup: kF \otimes k^Q \to k^Q$ and the invertible map $\sigma: F \times F \to (k^Q)^{\times}$ defined, respectively, by
\begin{equation} (x \rightharpoonup f) (p) = f(p \fiz x), \qquad x\in F, \, f \in k^Q, \, p \in Q; \end{equation}
\begin{equation}\label{def-sigma} \sigma_p(x, y) = \omega(p, x, y), \qquad x, y\in F, \,  p \in Q; \end{equation}
where $\sigma(x, y) = \sum_{p \in Q} \sigma_p(x, y) \delta_p$, for $x, y \in F$.
The normalized 3-cocycle condition \eqref{cc-omega} and the normalization assumption \eqref{norm-om} imply the following normalized 2-cocycle condition for $\sigma$:
\begin{align}\label{cc-sigma}& \sigma_{p \fiz x}(y, z)\sigma_p(x, yz)  = \sigma_p(xy, z)\sigma_p(x, y), \\ & \sigma_e(x, y)  = \sigma_p(e, y) = \sigma_p(x, e) = 1,
\end{align}
for all $x, y, z \in F$, $p \in Q$.
Thus $A$ is an associative algebra with unit element $\sum_{p \in Q}\delta_p \otimes 1$. For $f \in k^Q$, $x \in F$, the element $f\otimes x \in A$ will be denoted by $fx$. Hence, for all $x, y \in F$, $p, q \in Q$, we have
\begin{equation}\label{mult}(\delta_px) . (\delta_qy) :=  \delta_{p \fiz x, q} \sigma_p(x, y)  \, \delta_p xy. \end{equation}

\bigbreak
Consider the (non associative) crossed product  coalgebra structure on $A$ corresponding to the action $\fde$ and the invertible normalized map $\tau: Q \times Q \to (k^F)^{\times}$, given by
\begin{equation} \tau_x(p, q) = \omega(p, q, x), \qquad x\in F, \,  p, q \in Q; \end{equation}
where as before $\tau(p, q) = \sum_{x \in F} \tau_x(p, q) \delta_x$, for $p, q \in Q$. Using again the normalized 3-cocycle condition \eqref{cc-omega} and the normalization assumption \eqref{norm-om} on $\omega$, we find that $\tau$ satisfies the following 'twisted' normalized 2-cocycle condition:
\begin{align*}\label{cc-tau} \tau_x(p \fiz \theta(q, t), q . t) & \, \tau_x(q, t) \, \omega(p, q, t) \,
\sigma_p(\theta(q, t), q.t \fde x) \\ & \times \sigma_p(q \fde (t \fde x), \theta(q \fiz (t \fde x), t \fiz x))^{-1} \\
& = \tau_x(p . q, t) \, \tau_{t \fde x}(p, q) \, \omega(p \fiz (q \fde (t \fde x)), q \fiz (t \fde x), t \fiz x),
\end{align*}
\begin{equation*}\tau_e(p, q)  = \tau_x(e, q) = \tau_x(p, e) = 1,\end{equation*}
for all $x \in F$, $p, q, t \in Q$.

Explicitly, we have
\begin{equation}\label{comult}\Delta(\delta_px) := \sum_{s.t = p} \tau_x(s, t) \, \delta_s (t \fde x) \otimes \delta_tx, \qquad p \in Q, \, x\in F. \end{equation}
The counit for this coalgebra is given by $\epsilon \otimes \epsilon$.

\bigbreak
Both structures are related by the following theorem.

\begin{theorem} These algebra and coalgebra structures combine into a quasi-bialgebra structure on $A^{\op}$ with associator $\Phi \in A^{\otimes 3}$ given by
\begin{equation}\Phi = \sum_{p, q, r \in Q} \omega(p, q, r) \quad \delta_p\theta(q, r) \otimes \delta_q \otimes \delta_r. \end{equation}
There is a monoidal equivalence $\Rep (A^{\op}, \Phi) \sim \C(G, F, \omega, 1)$.
\end{theorem}

Note that $\Phi$ is invertible, with inverse $\Phi^{-1}$ given by the formula
\begin{align*}\Phi^{-1} = \sum_{p, q, r \in Q}  \omega(p, q, r)^{-1}  \sigma_p(\theta(q, r), & \theta(q, r)^{-1})^{-1} \\
& \delta_{p \fiz \theta(q, r)}\theta(q, r)^{-1} \otimes \delta_q \otimes \delta_r. \end{align*}

\pf Our definitions are dual to the ones given in Definition and Lemma 3.4.2 and Theorem and Definition 3.4.5 of \cite{sb}; note that, with the conventions of \cite{sb}, $\Phi$ is replaced by $\Phi^{-1}$ in condition \eqref{f4}.
Therefore, $A^{\op}$ is a quasi-bialgebra.
The monoidal equivalence $\Rep (A^{\op}, \Phi) \sim \C(G, F, \omega, 1)$ follows from \cite[Corollary 3.4.4]{sb}. \epf

\bigbreak
Since $\C(G, \omega, F, 1)$ is a rigid tensor category, it follows from \cite{sb-rig} that $A^{\op}$ is a quasi-Hopf algebra. We shall give the quasi-antipode in the next subsection.

Note that by Remark \ref{alfa=1} \emph{every} group theoretical category is equivalent to one of the form $\C(G, \omega, F, 1)$, for suitable $G, F$ and $\omega$, where $\omega$ satisfies \eqref{norm-om}, \eqref{norm-om-q}, in view of Proposition \ref{ln-om}. This gives us the following theorem.

\begin{theorem} Let $(H, \phi)$ be a finite dimensional quasi-Hopf algebra. Then $(H, \phi)$ is group theoretical if and only if it is gauge equivalent to a quasi-Hopf algebra of the form $(A^{\op}, \Phi)$, associated to suitable data $G$, $F$, $Q$ and $\omega$ satisfying \eqref{norm-om} and \eqref{norm-om-q}. \qed \end{theorem}

\bigbreak
We shall use the symbol $\circ$ to denote the multiplication in $A^{\op}$; so that $a \circ b = b . a$, for all $a, b \in A^{\op}$.

\bigbreak
\begin{obs} Using the properties listed in Lemma \ref{prop-beggs} and the normalization conditions \eqref{norm-om} and \eqref{norm-om-q}, it is not difficult to check that $A^{\op}$ is a quasi-bialgebra. For instance, $\Delta: A \otimes A \to A$ is an algebra map because of Lemma \ref{prop-beggs}-(ii), (iv) and the following relationship between $\sigma$ and $\tau$:
\begin{equation}\label{kac}\sigma_{t.s}(x, y) \tau_{xy}(t, s)  =
 \tau_x(t, s) \, \tau_y(t \fiz (s \fde x), s \fiz x) \,
\sigma_{t}(s \fde x, (s \fiz x) \fde y) \, \sigma_{s}(x, y), \end{equation}
for all $s, t \in Q$, $x, y \in F$,  which is a consequence of \eqref{norm-om-q} and \eqref{cc-omega}. Compare with \cite[Proposition 4.7]{ma-ext}.
\end{obs}

\begin{obs} Identify $\sigma$ and $\tau$, respectively, with maps
$$\sigma: Q \times F \times F \to k^{\times}, \qquad \tau: Q \times Q \times F \to k^{\times}.$$
Then the tuple $$(\Delta_G, 1, 1, .: Q \times Q \to Q, \fiz, \fde, \theta, \omega\vert_{Q \times Q \times Q}, \tau, \sigma),$$
constitutes the \emph{skeleton} of $(kG, \omega)$ according to \cite[Definition 4.1.1]{sb}.  \end{obs}

\subsection{} We give in this subsection the construction  of a quasi-antipode for $A^{\op}$.

We shall need the existence of inverses for the (nonassociative) multiplication in $Q$. This is guaranteed by the next lemma.

\begin{lema} The set $Q$ has well-defined left and right inverses with respect to the multiplication $.$; that is, for every $p \in Q$ there exist unique $p^L, p^R \in Q$ such that $p^L . p = e = p . p^R$. \end{lema}

Note that, by definition, we have
\begin{equation}pp^R = \theta(p, p^R), \qquad \text{and} \qquad p^Lp = \theta(p^L, p). \end{equation}

\pf As to left inverses, the lemma is contained in \cite[Proposition 2.3]{beggs}.  To prove the statement concerning right inverses, we shall use the assumption that $Q$ is also a set of representatives of the left cosets of $F$ in $G$.

Let $p \in Q$. By exactness of the factorization $G = QF$, there exist unique $s \in Q$, $x \in F$, such that $p^{-1} = sx$. Then we have
\begin{equation*} e = p p^{-1} = p s x = \theta(p, s) (p . s) x;
\end{equation*}
thus
\begin{equation*} \theta(p, s)^{-1} = (p . s) x.
\end{equation*}
Because $p.s \in Q$ and $\theta(p, s)^{-1}, x \in F$, the exactness of the factorization $G = QF$ implies that $p.s = e$.

We now show the uniqueness of such $s$, which gives the stament with $s = p^R$.  Suppose that $s' \in Q$ is such that $p.s' = e$. Then $ps', ps \in F$, and therefore also $(s')^{-1}s = (ps')^{-1} ps \in F$. This implies that $s' = s$, whence the uniqueness.  \epf

For later use, we give in the next lemma some of the relations between $(\quad)^L$, $(\quad)^R$ and the actions $\fde$, $\fiz$. The content of the lemma is part of \cite[Section 4]{beggs}.

\begin{lema}\label{pp-inv} The following relations hold, for all $p \in Q$:

(i) $p^{-1} = \theta(p^L, p)^{-1}p^L = p^R \theta(p, p^R)^{-1}$;

(ii) $p^L \fiz \theta(p, p^R) = p^R$ and $p^L \fde \theta(p, p^R) = \theta(p^L, p)$;

(iii) $p^{LL} = p \fiz \theta(p^L, p)^{-1}$;

(iv) $(p \fiz x)^L = p^L \fiz (p \fde x)$. \end{lema}

\pf The proof follows from the definitions and Lemma \ref{prop-beggs}.\epf

For notational convenience, we shall consider in the sequel the map $\succ$ introduced in the following definition.
Its main  properties are listed in the next lemma.

\begin{definition}\label{def-suc} The map $\succ: Q \times F \to F$ is defined as follows: $$p \succ x = \theta(p, p^R)^{-1} \, (p \fde x) \, \theta(p\fiz x, (p \fiz x)^R),$$
for all $p \in Q$, $x \in F$.
\end{definition}

\begin{lema}\label{pp-suc} Let $p \in Q$, $x, y \in F$. Then we have

(i) $p \succ (p^R \fde x) = x = p^R \fde (p \succ x)$;

(ii) $(p \fiz x)^R \fiz (p \succ x)^{-1} = p^R$;

(iii) $p \succ (xy) = (p \succ x) ((p \fiz x) \succ y)$;

(iv) $p \succ \theta(p^L, p)^{-1} = \theta(p, p^R)^{-1}$.
\end{lema}

\noindent In particular, for all $p \in Q$ the map $p \succ \underline{\quad}: F \to F$ is bijective, with inverse being $p^R \fde \underline{\quad}: F \to F$.

\pf We shall prove part (iv), the proof of (i)--(iii) being
straightforward. By Lemma \ref{pp-inv} (iii), we have $p
\theta(p^L, p)^{-1} = (p \fde \theta(p^L, p)^{-1}) p^{LL}$, and on
the other hand, $p \theta(p^L, p)^{-1} = (p^L)^{-1} =
\theta(p^{LL}, p^L)^{-1} p^{LL}$.

By exactness of the
factorization $G = FQ$, we get $\theta(p^{LL}, p^L)^{-1} = p \fde
\theta(p^L, p)^{-1}$. Now, by definition, $$p \succ \theta(p^L,
p)^{-1} = \theta(p, p^R)^{-1} \, (p \fde \theta(p^L, p)^{-1}) \,
\theta(p^{LL}, p^L) = \theta(p, p^R)^{-1},$$ by the above. This
proves (iv). \epf

\begin{theorem}\label{q-antip} There is a quasi-Hopf algebra structure on $A^{\op}$, with quasi-antipode $\mathcal S: A^{\op} \to A^{\op}$ given by
\begin{equation}\label{antip}\mathcal S(\delta_px) = \tau_{p \succ x}(p, p^R)^{-1} \, \sigma_{p^R}(p\succ x, (p\succ x)^{-1})^{-1}
\, \delta_{(p\fiz x)^R} (p\succ x)^{-1}.
\end{equation}
We have $\alpha = 1$ and $$\beta = \sum_{q \in Q} \omega(q^{-1}, q, q^{-1}) \, \delta_q \theta(q^L, q)^{-1} = \sum_{q \in Q} \omega(q, q^{-1}, q)^{-1} \, \delta_q \theta(q^L, q)^{-1}.$$ \end{theorem}

Compare with the formulas given in \cite{zhu}, \cite{beggs} for the case where $\omega = 1$.

\pf We shall freely use the relations in Lemma \ref{eq-om} and the normalization conditions \eqref{norm-om} and \eqref{norm-om-q}.
Using relation \eqref{kac} and Lemma \ref{pp-suc}, it is straightforward to see that $\mathcal S$ is an anti-algebra map.
The injectivity of $\mathcal S$ follows from the injectivity of the map $p \succ \underline{\quad}: F \to F$ and relation (ii) in Lemma \ref{pp-suc}. Therefore $\mathcal S: A^{\op} \to A^{\op}$ is an algebra anti-automorphism.

\bigbreak
We now check condition \eqref{ant-ab}.
Let $p \in Q$, $x \in F$, and let $X = \delta_px \in A^{\op}$. We have
\begin{flalign*}& \mathcal S(X_1) \circ \alpha \circ X_2  = \mathcal S(X_1) \circ  X_2 = X_2 . \mathcal S(X_1) \\
& =  \sum_{s.t = p} \tau_x(s, t) \, \delta_tx \, .  \mathcal S(\delta_s (t \fde x)) \\
& = \sum_{s.t = p} \tau_x(s, t) \, \tau_{s \succ (t \fde x)}(s, s^R)^{-1} \, \sigma_{s^R}(s\succ (t \fde x), (s\succ (t \fde x))^{-1})^{-1}
\\ & \times \delta_tx \, . \delta_{(s\fiz (t \fde x))^R} (s\succ (t \fde x))^{-1} \\
& = \sum_{s.t = p} \tau_x(s, t) \, \tau_{s \succ (t \fde x)}(s, s^R)^{-1} \, \sigma_{s^R}(s\succ (t \fde x), (s\succ (t \fde x))^{-1})^{-1}
\\ & \times \delta_{t \fiz x, (s\fiz (t \fde x))^R} \quad \sigma_t(x, (s \succ (t \fde x))^{-1}) \quad \delta_t \, x (s\succ (t \fde x))^{-1}. \end{flalign*}

By Lemma \ref{pp-inv} (iv), we have $\delta_{t \fiz x, (s\fiz (t \fde x))^R} = \delta_{s, t^L}$. Hence, using property (i) in Lemma \ref{pp-suc}, the last expression equals
\begin{flalign*}\delta_{p, e} & \sum_{s \in Q} \tau_x(s, s^R) \, \tau_{s \succ (s^R \fde x)}(s, s^R)^{-1} \, \sigma_{s^R}(s\succ (s^R \fde x), (s\succ (s^R \fde x))^{-1})^{-1}
\\ & \times \sigma_{s^R}(x, (s \succ (s^R \fde x))^{-1}) \quad \delta_{s^R} \, x (s\succ (s^R \fde x))^{-1} \\
& = \delta_{p, e} \sum_{s \in Q} \delta_{s^R} = \delta_{p, e} 1 =
\delta_{p, e} \alpha. \end{flalign*} This proves the right hand
side identity in \eqref{ant-ab}. We now compute
\begin{flalign*}& X_1 \circ \beta \circ \mathcal S(X_2) =  \mathcal S(X_2) . \beta . X_1 =  \sum_{s.t = p} \tau_x(s, t) \, \mathcal S(\delta_tx) \, . \beta \, . \delta_s (t \fde x) \\
& =  \sum_{s.t = p}\sum_{q} \omega(q, q^{-1}, q)^{-1} \, \tau_x(s, t) \, \mathcal S(\delta_tx) \, . \delta_q \theta(q^L, q)^{-1} \, . \delta_s (t \fde x) \\
& =  \sum_{s.t = p}\sum_{q} \omega(q, q^{-1}, q)^{-1} \,  \tau_x(s, t) \, \tau_{t \succ x}(t, t^R)^{-1} \, \sigma_{t^R}(t\succ x, (t\succ x)^{-1})^{-1} \\
& \times \delta_{(t\fiz x)^R} (t\succ x)^{-1} \, .
\delta_q \theta(q^L, q)^{-1} \, . \delta_s (t \fde x) \\
& =  \sum_{s.t = p}\sum_{q} \omega(q, q^{-1}, q)^{-1} \,  \tau_x(s, t) \, \tau_{t \succ x}(t, t^R)^{-1} \, \sigma_{t^R}(t\succ x, (t\succ x)^{-1})^{-1} \\
& \times \delta_{(t\fiz x)^R \fiz (t\succ x)^{-1}, q} \,
\sigma_{(t\fiz x)^R}((t\succ x)^{-1}, \theta(q^L, q)^{-1}) \\
& \times \delta_{(t\fiz x)^R} \, (t\succ x)^{-1} \theta(q^L, q)^{-1} \, \delta_s (t \fde x). \end{flalign*}
By Lemma \ref{pp-suc} (ii), this equals
\begin{flalign*} \sum_{s.t = p} & \omega(t^R, (t^R)^{-1}, t^R)^{-1} \,  \tau_x(s, t) \, \tau_{t \succ x}(t, t^R)^{-1} \, \sigma_{t^R}(t\succ x, (t\succ x)^{-1})^{-1} \\
&  \sigma \sigma_{(t\fiz x)^R}((t\succ x)^{-1}, \theta(t, t^R)^{-1}) \,
 \delta_{(t\fiz x)^R} \, (t\succ x)^{-1} \theta(t, t^R)^{-1} \, \delta_s (t \fde x)\\
 =  \sum_{s.t = p} & \omega(t^R, (t^R)^{-1}, t^R)^{-1} \,  \tau_x(s, t) \, \tau_{t \succ x}(t, t^R)^{-1} \, \sigma_{t^R}(t\succ x, (t\succ x)^{-1})^{-1} \\
&  \times \sigma_{(t\fiz x)^R}((t\succ x)^{-1}, \theta(t, t^R)^{-1}) \,
\sigma_{(t\fiz x)^R}((t\succ x)^{-1} \theta(t, t^R)^{-1}, t \fde x) \\
& \times \delta_{(t\fiz x)^R \fiz (t\succ x)^{-1} \theta(t, t^R)^{-1}, s} \,
 \delta_{(t\fiz x)^R} \,
(t\succ x)^{-1} \theta(t, t^R)^{-1} (t \fde x) \\
= \delta_{p, e} &  \sum_{t} \omega(t^R, (t^R)^{-1}, t^R)^{-1} \,  \tau_x(t^L, t) \, \tau_{t \succ x}(t, t^R)^{-1} \\ & \times \sigma_{t^R}(t\succ x, (t\succ x)^{-1})^{-1} \,
  \sigma_{(t\fiz x)^R}((t\succ x)^{-1}, \theta(t, t^R)^{-1}) \\
& \times \sigma_{(t\fiz x)^R}((t\succ x)^{-1} \theta(t, t^R)^{-1}, t \fde x) \quad  \delta_{(t\fiz x)^R} \, \theta(t \fiz x, (t \fiz x)^R)^{-1}; \end{flalign*}
the last equality by Lemma \ref{pp-suc} (ii) and Definition \ref{def-suc}.

Using the cocycle condition \eqref{cc-sigma} and Lemma \ref{pp-inv}, we find
\begin{flalign*}& \sigma_{t^R}(t\succ x, (t\succ x)^{-1})^{-1} \,
  \sigma_{(t\fiz x)^R}((t\succ x)^{-1}, \theta(t, t^R)^{-1}) \\
& \times \sigma_{(t\fiz x)^R}((t\succ x)^{-1} \theta(t, t^R)^{-1}, t \fde x) \\
& = \sigma_{t^R}(\theta(t, t^R)^{-1}, t\fde x) \,
  \sigma_{t^R}(t\succ x, \theta(t \fiz x, (t \fiz x)^R)^{-1})^{-1}\\
& = \omega(t^R, \theta(t, t^R)^{-1}, t\fde x) \\
& \times \omega(t^R, \theta(t, t^R)^{-1} (t\fde x) \theta(t \fiz x, (t \fiz x)^R), \theta(t \fiz x, (t \fiz x)^R)^{-1})^{-1} \\
& = \omega(t^R, \theta(t, t^R)^{-1}, t\fde x) \,
\omega((t \fiz x)^L, \theta(t \fiz x, (t \fiz x)^R),  \theta(t \fiz x, (t \fiz x)^R)^{-1})^{-1} \\
& \times \omega(t^R, \theta(t, t^R)^{-1}(t\fde x), \theta(t \fiz x, (t \fiz x)^R))\\
& = \omega((t \fiz x)^L, \theta(t \fiz x, (t \fiz x)^R),  \theta(t \fiz x, (t \fiz x)^R)^{-1})^{-1} \\
& \times \omega(t^L, t\fde x, \theta(t \fiz x, (t \fiz x)^R)) \,
\omega(t^R, \theta(t, t^R)^{-1}, (t\fde x) \theta(t \fiz x, (t \fiz x)^R))\\
& = \omega((t \fiz x)^L, \theta(t \fiz x, (t \fiz x)^R),  \theta(t \fiz x, (t \fiz x)^R)^{-1})^{-1} \\
& \times \omega((t \fiz x)^L, t \fiz x, (t \fiz x)^R))^{-1} \, \omega(t^L, t, x)^{-1} \, \omega(t^L, t, x(t \fiz x)^R) \\
& \times \omega(t^R, \theta(t, t^R)^{-1}, (t\fde x) \theta(t \fiz
x, (t \fiz x)^R)); \end{flalign*} where we have used the
definition of $\sigma$ \eqref{def-sigma}, the cocycle condition on
$\omega$, and Lemma \ref{eq-om} (iii).

On the other hand,
\begin{flalign*}& \omega((t \fiz x)^L, \theta(t \fiz x, (t \fiz x)^R),  \theta(t \fiz x, (t \fiz x)^R)^{-1})^{-1} \,
 \omega((t \fiz x)^L, t \fiz x, (t \fiz x)^R))^{-1}\\
& = \omega((t \fiz x)^R, \theta(t \fiz x, (t \fiz x)^R)^{-1},  (t \fiz x) (t \fiz x)^R)^{-1} \,
 \omega((t \fiz x)^L, t \fiz x, (t \fiz x)^R))^{-1}\\
&  = \omega((t \fiz x)^{-1}, t \fiz x, (t \fiz x)^R) \,
\omega((t \fiz x)^R, ((t \fiz x)^R)^{-1}, (t \fiz x)^R)^{-1} \\
& \times \omega((t \fiz x)^L, t \fiz x, (t \fiz x)^R))^{-1}\\
& = \omega((t \fiz x)^R, ((t \fiz x)^R)^{-1}, (t \fiz x)^R)^{-1}.
\end{flalign*}
Therefore
\begin{flalign*}& \omega(t^R, (t^R)^{-1}, t^R)^{-1} \,  \tau_x(t^L, t) \, \tau_{t \succ x}(t, t^R)^{-1} \, \sigma_{t^R}(t\succ x, (t\succ x)^{-1})^{-1} \\
&  \times \sigma_{(t\fiz x)^R}((t\succ x)^{-1}, \theta(t, t^R)^{-1}) \,
 \sigma_{(t\fiz x)^R}((t\succ x)^{-1} \theta(t, t^R)^{-1}, t \fde x) \\
& = \omega(t^R, (t^R)^{-1}, t^R)^{-1} \, \omega(t, t^R, t \succ x)^{-1} \, \omega((t \fiz x)^R, ((t \fiz x)^R)^{-1}, (t \fiz x)^R)^{-1} \\
&  \times  \omega(t^{-1}, t, x(t \fiz x)^R)  \,  \omega(t^R, \theta(t, t^R)^{-1}, (t\fde x) \theta(t \fiz x, (t \fiz x)^R)) \\
& = \omega(t^R, (t^R)^{-1}, t^R)^{-1} \, \omega(t, t^R \theta(t, t^R)^{-1}, (t \fde x) \theta(t \fiz x, (t \fiz x)^R))^{-1} \\
&  \times \omega(t, t^R, \theta(t, t^R)^{-1})^{-1} \, \omega((t \fiz x)^R, ((t \fiz x)^R)^{-1}, (t \fiz x)^R)^{-1} \,  \omega(t^L, t, x(t \fiz x)^R) \\
& = \omega(t^R, (t^R)^{-1}, t^R)^{-1} \, \omega(t, t^{-1}, t x (t \fiz x)^R)^{-1} \\
& \times  \omega(t, t^R, \theta(t, t^R)^{-1})^{-1} \, \omega((t \fiz x)^R, ((t \fiz x)^R)^{-1}, (t \fiz x)^R)^{-1} \,  \omega(t^{-1}, t, x(t \fiz x)^R) \\
& = \omega(t^R, (t^R)^{-1}, t^R)^{-1} \, \omega(t, t^{-1}, t)^{-1}  \, \omega(t^{-1}, t, x(t \fiz x)^R)^{-1}\\
&  \times \omega(t, t^R, \theta(t, t^R)^{-1})^{-1} \, \omega((t \fiz x)^R, ((t \fiz x)^R)^{-1}, (t \fiz x)^R)^{-1} \,  \omega(t^L, t, x(t \fiz x)^R) \\
& = \omega(t^R, (t^R)^{-1}, t^R)^{-1} \, \omega(t, t^R \theta(t, t^R)^{-1}, t)^{-1}   \,  \omega(t, t^R, \theta(t, t^R)^{-1})^{-1} \\
&  \times  \omega((t \fiz x)^R, ((t \fiz x)^R)^{-1}, (t \fiz x)^R)^{-1} \\
& = (\omega(t^R, (t^R)^{-1}, t^R) \omega(t, t^R, \theta(t,
t^R)^{-1} t)
 \omega((t \fiz x)^R, ((t \fiz x)^R)^{-1}, (t \fiz x)^R))^{-1} \\
& = \omega(t^R, (t^R)^{-1}, t^R)^{-1} \, \omega((t^R)^{-1}, t^R,
(t^R)^{-1})^{-1} \\ & \omega((t \fiz x)^R, ((t \fiz x)^R)^{-1}, (t
\fiz x)^R)^{-1} \\ & = \omega((t \fiz x)^R, ((t \fiz x)^R)^{-1},
(t \fiz x)^R)^{-1}.  \end{flalign*} Hence we get
\begin{flalign*} X_1 \circ \beta \circ \mathcal S(X_2)  = & \delta_{p, e}   \sum_{t} \omega((t \fiz x)^R, ((t \fiz x)^R)^{-1}, (t \fiz x)^R)^{-1} \\ &  \times  \delta_{(t\fiz x)^R} \, \theta(t \fiz x, (t \fiz x)^R)^{-1} = \delta_{p, e} \beta, \end{flalign*}
which gives the right hand side identity in \eqref{ant-ab}.

\bigbreak
Finally, the proof of conditions \eqref{largas} is straightforward, using the properties listed in Lemma \ref{pp-inv} and the cocycle conditions. This finishes the proof of the theorem.
\epf

\section{Frobenius-Schur indicators for $\C(G, \omega, F, 1)$}\label{fs}

Let $G$ be a finite group, $F \subseteq G$  a subgroup, and $\omega: G \times G \times G \to k^{\times}$ a 3-cocycle subject to the normalization conditions \eqref{norm-om} and \eqref{norm-om-q}. We keep the notation of the previous sections for the skeleton maps $\fde$, $\fiz$, $\theta$, $\sigma$ and $\tau$.

We have a $k$-linear monoidal equivalence $\Rep (A^{\op}, \Phi) \sim \C(G, F, \omega, 1)$, where $(A^{\op}, \Phi)$ is the quasi-Hopf algebra attached to the data $(G, F, \omega)$ in Section \ref{qha}.
Our aim in this section is to give an explicit description of the canonical central  element $\nu_{A^{\op}} \in A^{\op}$ and then of the Frobenius-Schur indicators for the quasi-Hopf algebra $A^{\op}$. It follows from gauge invariance of the Frobenius-Schur indicators that these depend only on the fusion category $\C(G, F, \omega, 1)$.

\subsection{}
Let $\Lambda_0 : = \vert F \vert^{-1} \sum_{x \in F}x \in kF$ be the normalized integral.
The normalized two sided integral $\Lambda \in A^{\op}$ has the following form:
\begin{equation}\label{integral}\Lambda = \delta_e \Lambda_0 = \vert F \vert^{-1} \sum_{x \in F} \delta_e x. \end{equation}

\begin{prop}\label{beta-inv} The element $\beta$ is invertible with inverse
$$\beta^{-1} = \sum_{p \in Q} \omega(p^L, p, p^R) \, \delta_{p^L} \theta(p, p^R).$$
We have in addition $\mathcal S(\beta) = \beta^{-1}$.
\end{prop}

\pf It is not difficult to check that the expression
$$\sum_{p \in Q} \omega(p^R, (p^R)^{-1}, p^R) \, \sigma_{p^R}(\theta(p, p^R)^{-1}, \theta(p, p^R))^{-1}
\, \delta_{p^L} \theta(p, p^R).$$
defines an inverse of $\beta$. We claim that
\begin{equation}\label{id-om}\omega(p^R, (p^R)^{-1}, p^R) \, \sigma_{p^R}(\theta(p, p^R)^{-1}, \theta(p, p^R))^{-1} = \omega(p^L, p, p^R),\end{equation}
for all $p \in Q$. This will imply the claimed expression for
$\beta^{-1}$.

Letting $q = p^R$, equation \eqref{eq-om} is equivalent to the following:
\begin{equation}\label{id-om-2}\omega(q, q^{-1}, q) \, \sigma_q(\theta(q^L, q)^{-1}, \theta(q^L, q))^{-1}=  \omega(q \fiz \theta(q^L, q)^{-1}, q^L, q). \end{equation}

To establish equation \eqref{id-om-2}, we note that $\theta(q^L, q) = q^L q \in F$, for all $q \in Q$. Then, applying the cocycle and normalization conditions on $\omega$, we get
\begin{align*} \omega(q \fiz \theta(q^L, q)^{-1}, q^L, q) & =
\omega(q, q^{-1}, q) \, \omega(q, \theta(q^L, q)^{-1}, q^Lq)^{-1} \\
& = \omega(q, q^{-1}, q) \, \omega(q, \theta(q^L, q)^{-1}, \theta(q^L, q))^{-1}, \end{align*}
which is  the claimed identity.

\bigbreak Using Lemma \ref{pp-suc}(iv), we get
\begin{align*}
\mathcal S(\delta_p \theta(p^L, p)^{-1})  = & \tau_{\theta(p,
p^R)^{-1}} (p, p^R)^{-1} \, \sigma_{p^R}(\theta(p, p^R)^{-1},
\theta(p, p^R))^{-1} \\ & \times \delta_{p^L} \, \theta(p, p^R).
\end{align*}
We now compute
\begin{align*}
\tau_{\theta(p, p^R)^{-1}} & (p, p^R)  = \omega(p, p^R,
(pp^R)^{-1}) \\
& = \omega(pp^R, (p^R)^{-1}, p^{-1})^{-1} \, \omega(p, p^R,
(p^R)^{-1})
\, \omega(p^R, (p^R)^{-1}, p^{-1}) \\
& = \omega(p, p^R, (p^R)^{-1}) \, \omega(p^R, (p^R)^{-1}, p^{-1}),
\end{align*} because $pp^R \in F$. Similarly,
\begin{align*}
\sigma_{p^R}(\theta(p, & p^R)^{-1}, \theta(p, p^R))  = \omega(p^R,
(p^R)^{-1}p^{-1}, pp^R) \\
& = \omega(p^R, (p^R)^{-1}, p^{-1})^{-1} \, \omega((p^R)^{-1},
p^{-1}, pp^R)^{-1}
\, \omega(p^R, (p^R)^{-1}, p^R) \\
& = \omega(p^R, (p^R)^{-1}, p^R) \, \omega(p, p^{-1},
(pp^R))^{-1}.
\end{align*} Hence
\begin{align*}
\tau_{\theta(p, p^R)^{-1}} (p, p^R) & \, \sigma_{p^R}(\theta(p,
p^R)^{-1}, \theta(p, p^R)) = \omega(p, p^{-1}, p)^{-1} \,
\omega(p^{-1}, p, p^R)^{-1} \\
& = \omega(p, p^{-1}, p)^{-1} \, \omega(p^L, p, p^R)^{-1}.
\end{align*}
Thus,
\begin{align*}
\mathcal S(\beta) & = \sum_p \omega(p^{-1}, p, p^{-1}) \, \mathcal
S(\delta_p \theta(p^L, p)^{-1}) \\
& = \sum_p\omega(p^L, p, p^R) \, \delta_{p^L} \theta(p, p^R) \\
& = \beta^{-1}.
\end{align*}
This finishes the proof of the proposition. \epf

\begin{theorem}\label{cc-el} The canonical central  element $\nu_{A^{\op}}$ is given by the formula
\begin{equation}\label{nu-1}\nu_{A^{\op}}  =
\vert F \vert^{-1}   \sum_{(q \fiz x).q = e} \omega(xq, xq, xq) \quad \delta_{q} (xq)^2.  \end{equation}
We have also
\begin{flalign*}\nu_{A^{\op}}  =
\vert F \vert^{-1}  & \sum_{(q \fiz x).q = e} \tau_x(q \fiz x, q) \, \sigma_q(x, q \fde x) \, \omega((q \fiz x)^L, q \fiz x, q) & \\
& \times \sigma_q(x(q \fde x), (q \fiz x)q) \quad \delta_{q} (xq)^2. & \end{flalign*}
\end{theorem}

\pf Since the element $\beta$ corresponding to the quasi-antipode of $A^{\op}$ is invertible and $\alpha = 1$, we have
$$\nu_{A^{\op}} = (\Lambda_1 \circ \Lambda_2) \circ \beta^{-1} =  \beta^{-1} \circ (\Lambda_1 \circ \Lambda_2),$$
and thus
$$\nu_{A^{\op}} = \beta^{-1} \circ (\Lambda_1 \circ \Lambda_2)   =  (\Lambda_1 \circ \Lambda_2) . \beta^{-1}.$$
Using formula \eqref{integral} for the integral $\Lambda$, we find
\begin{equation}\label{delta-lambda}\Delta(\Lambda) = \Lambda_1 \otimes \Lambda_2 =
\vert F \vert^{-1} \sum_{x \in F} \sum_{q \in Q} \tau_x(q^L, q) \quad
\delta_{q^L} (q \fde x) \otimes \delta_qx,\end{equation}
thus
\begin{equation}\Lambda_1 \circ \Lambda_2 = \Lambda_2 . \Lambda_1 =
\vert F \vert^{-1}  \sum_{(q \fiz x).q = e} \tau_x(q \fiz x, q) \, \sigma_q(x, q \fde x)  \quad \delta_q x(q \fde x). \end{equation}
From these, we  compute $\beta^{-1} \circ (\Lambda_1 \circ \Lambda_2)$ and get the second expression for $\nu_{A^{\op}}$.
Now, for all $x \in F$ and $q \in Q$ such that $(q \fiz x).q = e$, we have
\begin{align*}
& \sigma_q(x, q \fde x) \, \sigma_q(x(q \fde x), (q \fiz x)q) =
\sigma_{q\fiz x}(q \fde x, (q \fiz x)q) \, \sigma_q(x, qxq) \\
& = \omega(q\fiz x, q \fde x, (q \fiz x)q) \, \omega(q, x, qxq) \\
& = \omega(q\fiz x, q \fde x, (q \fiz x)q) \, \omega(q\fiz x, q, xq)^{-1} \, \omega(q, xq, xq) \\
& = \omega(q\fiz x, q \fde x, (q \fiz x)q) \, \omega(q\fiz x, q, x)^{-1} \,  \omega(qx, qx, q)^{-1} \, \omega(q, xq, xq).\end{align*}
The last equality because
\begin{align*} \omega(q\fiz x, q, xq) &  =  \omega(q\fiz x, q, x) \,   \omega(q\fiz x, qx, q) \\
& =  \omega(q\fiz x, q, x) \,  \omega(qx, qx, q).\end{align*}

On the other hand,
\begin{align*}
\omega((q \fiz x)^L, q \fiz x, q) & = \omega(q^L \fiz (q \fde x), q \fiz x, q) \\ & = \omega(q^L, qx, q) \, \omega(q^L, q \fde x, (q\fiz x)q)^{-1} \\
& = \omega(q\fiz x, qx, q) \, \omega(q\fiz x, q \fde x, (q\fiz x)q)^{-1} \\ & = \omega(qx, qx, q) \, \omega(q\fiz x, q \fde x, (q\fiz x)q)^{-1}. \end{align*}
Therefore
\begin{align*} \tau_x(q \fiz x, q) \, \sigma_q(x, q \fde x) & \, \omega((q \fiz x)^L, q \fiz x, q) \,  \sigma_q(x(q \fde x), (q \fiz x)q)  \\ & = \omega(q, xq, xq) = \omega(xq, xq, xq).\end{align*}
This proves equation \eqref{nu-1} and finishes the proof of the proposition. \epf

\begin{obs}Computing instead $(\Lambda_1 \circ \Lambda_2)  \circ \beta^{-1}$, we get
\begin{flalign*}\label{nu-2}\nu_{A^{\op}}  =
\vert F \vert^{-1}  & \sum_{(q \fiz x).q = e} \tau_x(q \fiz x, q) \, \sigma_q(x, q \fde x) \, \omega((q \fiz x)^L, q \fiz x, q) & \\
& \omega((q \fiz x)^L, (q \fiz x)q, x(q \fde x)) \quad \delta_{(q \fiz x)^L} ((q \fiz x)(q \fde x))^2. & \end{flalign*}\end{obs}

As a consequence of Theorem \ref{cc-el}, we get the following expression for the Frobenius-Schur indicators.
After suitable normalization, this expression allows to compute the Frobenius-Schur indicators for every group theoretical category.

\begin{cor}\label{fs-chi} Suppose $\chi$ is an irreducible character of $A^{\op}$.  Then the Frobenius-Schur indicator of $\chi$ is given by
\begin{flalign*}\chi(\nu_{A^{\op}}) & =
\vert F \vert^{-1}   \sum_{(q \fiz x).q = e} \omega(xq, xq, xq)  \quad \chi(\delta_{q} (xq)^2) \\
& = \vert F \vert^{-1} \sum_{(q \fiz x).q = e} \tau_x(q \fiz x, q) \, \sigma_q(x, q \fde x) \, \omega((q \fiz x)^L, q \fiz x, q) & \\
& \qquad \qquad \sigma_q(x(q \fde x), (q \fiz x)q) \quad
\chi(\delta_{q} (xq)^2). & \qed \end{flalign*}
\end{cor}

\subsection{} In this subsection we aim to give an explicit description of the irreducible characters (and hence of the indicators) of
$\C(G, \omega, F, 1)$ in terms of the groups $G$ and $F$.

As an algebra $A = k^Q \#_{\sigma}kF$ is a crossed product. See Subsection \ref{str-qha}. Hence the irreducible left $A$-modules can be described using Clifford theory.

On the other hand, to every left $A$-module $V$ one can associate the left $A^{\op}$-module $V^*$, the action of $a \in A^{\op}$ being the transpose of the action of $a \in A$ on $V$. This gives a bijective correspondence between (irreducible) left $A$-modules $V$ and (irreducible) left $A^{\op}$-modules. Moreover, this bijection preserves characters: $\chi_{V^*} = \chi_V$, for all finite dimensional left $A$-module $V$.

\bigbreak
Let $F^p \subseteq F$ denote the isotropy subgroup of $p \in Q$. Then
the restriction of $\sigma_p$  defines a normalized  2-cocycle  $$\sigma_p: F^p \times F^p \to k^{\times}.$$
Let $k_{\sigma_p}F^p$ denote the corresponding twisted group algebra.

The space of isomorphism classes of irreducible $A$-modules can be parametrized by the modules $V_{p, W}$, where
\begin{equation}V_{p, W} = \Ind_{k^Q \#_{\sigma}kF^p} \, p \otimes W = A \otimes_{k^Q \#_{\sigma}kF^p} (p \otimes W),
\end{equation}
where $p$ runs over a set of representatives of the action of $F$ on $Q$, and $W$ runs over a system of representatives of isomorphism classes of irreducible left $k_{\sigma_p}F^p$-modules. See \cite[Section 3]{KMM}.

\bigbreak
There is a natural identification between $Q$ and the space $F\backslash G = \{ Fg: g \in G \}$ of left cosets of $F$ in $G$.
Under this identification, the action $Q \times F \to Q$ corresponds to the natural action of $F$ on $F\backslash G$ by right multiplication: $Fg . x = F(gx)$, $g \in G$, $x \in F$.

This gives in turn a natural identification between the space of orbits of the action  $Q \times F \to Q$ and the space $F\backslash G / F$ of double cosets of $F$ in $G$.
Moreover, the isotropy subgroup of an element $p \in Q$ is $F^p = F \cap p^{-1}Fp$. Hence we get

\begin{prop}\label{irr-mod} The set of isomorphism classes of irreducible $A^{\op}$-modules is parametrized by the modules $U_{p, W}$, where
\begin{equation}U_{p, W} = V_{p, W}^* = (\Ind_{k^Q \#_{\sigma}kF^p} \, p \otimes W)^*,
\end{equation}
where $p$ runs over a set of representatives of the double cosets of $F$ in $G$, $F^p = F \cap p^{-1}Fp$, and $W$ runs over a system of representatives of isomorphism classes of irreducible left $k_{\sigma_p}F^p$-modules.

\bigbreak
The character of the irreducible $A^{\op}$-module $U_{p, W}$ is given by the formula
\begin{equation}\chi_{p, W}(\delta_qz) = \sum_{y^{-1}zy \in F^p} \delta_{p, q \fiz y} \, \sigma_q(z, y)\, \sigma_q(y, y^{-1}zy)^{-1} \, \chi_W (y^{-1}zy),
\end{equation}
where the sum is over all $y$ running over a set of representatives of the right cosets of $F^p$ in $F$, and $\chi_W$ is the character of $W$.
\end{prop}

\bigbreak Observe that $\dim U_{p, W} = [F: F \cap p^{-1}Fp] \dim
W$. So the Proposition immediately implies that the dimensions of
the irreducible modules of a group-theoretical quasi-Hopf algebra
divide its dimension, \textit{i.e.}, that Kaplanksy's conjecture
holds in this case.

\bigbreak
\begin{proof} We only need to prove the formula for the character.
The character of $U_{p, W}$ coincides with the character of $V_{p,
W}$.  Let $Y$ be a set of representatives of the right cosets of
$F^p$ in $F$. A basis of $V_{p, W}$ is given by $y \otimes p
\otimes v$, where $(v)$ is a basis of $W$, and $y \in Y$.

For all $q \in Q$, $y, z \in F$, we have
\begin{align*}(\delta_qz) . y & = \sigma_q(z, y) \, \delta_q zy \\
& = \sigma_q(z, y) \, \delta_q y(y^{-1}zy) \\
& = \sigma_q(z, y) \, \sigma_q(y, y^{-1}zy)^{-1} \, y . (\delta_{q \fiz y} y^{-1}zy. \\
\end{align*}
Hence, the action of $\delta_qz$ on this basis is
\begin{align*}(\delta_qz) . y \otimes p \otimes v & = (\delta_qz) . y \otimes
p \otimes v \\ & = \sigma_q(z, y) \, \sigma_q(y, y^{-1}zy)^{-1} \,
y . (\delta_{q \fiz y} y^{-1}zy) \otimes p \otimes v.\end{align*}
Thus, in order to compute the trace of this action, we  only need
to consider those basis vectors $y \otimes p \otimes v$, for which
$y^{-1}zy \in F^p$; and for such  $y$, we have
\begin{align*}(\delta_qz) . y \otimes p \otimes v  & = \sigma_q(z,
y)\, \sigma_q(y, y^{-1}zy)^{-1} \, y (\delta_{q \fiz y} y^{-1}zy) \otimes p \otimes v \\
& = \delta_{p, q \fiz
y} \, \sigma_q(z, y)\, \sigma_q(y, y^{-1}zy)^{-1} \, y \otimes p
\otimes (y^{-1}zy) . v.\end{align*} This implies the desired
formula.
\end{proof}

\bigbreak
\begin{obs} The parametrization in Proposition \ref{irr-mod} allows to recover the statement in the Remark after Proposition 3.1 of \cite{ostrik}, for the category $\C(G, \omega, F, 1)$.
Indeed, the 2-cocycle $\psi^p(x, y) \in Z^2(F^p, k^{\times})$ considered in {\it loc. cit.} coincides in our notation with $\sigma_p(y^{-1}, x^{-1})$; and this is cohomologous to $\sigma_p(x, y)$ via $d(\gamma)$, where $\gamma(x) = \sigma_p(x^{-1}, x)$, $x \in F$.
\end{obs}

\section{Examples}\label{ejemplos}

In this section we discuss some special cases of the results in Sections \ref{qha}, \ref{fs}.

\subsection{Abelian extensions}
Suppose that $G = FQ$ is an \emph{exact factorization} of the group $G$; that is, $Q$ is a subgroup of $G$ and $(F, Q)$ is a \emph{matched pair} of finite groups with the actions $\fde : Q \times F \to F$, $\fiz : Q \times F \to Q$. We refer the reader to \cite{ma-ext, ma-ext2} for the main notions used here, and in particular for the study of the cohomology theory associated to the matched pair $(F, Q)$.

\bigbreak
Fix a representative $(\tau, \sigma)$
of a class in $\Opext (k^G, kF)$; that is, $\sigma: F \times F \to (k^Q)^{\times}$ and $\tau: Q \times Q \to (k^F)^{\times}$ are normalized 2-cocycles subject to compatibility conditions.

Consider the 3-cocycle $\omega: G \times G \times G \to k^{\times}$   given by
\begin{equation}\label{om-mp} \omega (\tau, \sigma) \left(xp, yq, zr \right) =
\tau_{z}(p \fiz y, q) \, \sigma_{p}(y,
q \fde z), \end{equation} for all $x, y, z \in F$, $p, q, r \in Q$.
The class of the cocycle $\omega = \omega (\tau, \sigma)$ is the image of the class of $(\tau, \sigma)$ in the Kac exact sequence \cite{sb, ma-ext2}.

It is not difficult to see that $\sigma$ and $\tau$ have the same meaning as in Subection \ref{str-qha}. Note that $\omega\vert_{Q \times Q \times Q} = 1$.

\bigbreak
There is a bicrossed product Hopf algebra $A : = k^G
\, {}^{\tau}\#_{\sigma}kF$ corresponding to this data. As is well-known, this correspondence gives a bijection between the equivalence classes of Hopf algebra extensions  $$1 \to k^Q \to A \to kF \to 1,$$ and the abelian group $\Opext (k^G, kF)$.
The Hopf algebra $A^{\op}$ coincides with the (quasi-)Hopf algebra corresponding to $G$, $F$ and $\omega$, as in Subsection \ref{str-qha}.

Applying Corollary \ref{fs-chi}, we find the following expression for the Frobenius-Schur indicators.

\begin{prop} Let $\chi$ be an irreducible character of $A^{\op}$. Then the Frobenius-Schur indicator of $\chi$ is given by
\begin{align*}\chi(\nu_{A^{\op}}) & = \vert F \vert^{-1} \sum_{q \fiz x = q^{-1}} \tau_x(q^{-1}, q) \, \sigma_q(x, q \fde x) \,
\quad \chi (\delta_{q} x(q\fde x)) \\
& = \vert F \vert^{-1} \sum_{q \fiz x = q^{-1}} \tau_x(q^{-1}, q) \, \sigma_q(x, qxq) \, \quad \chi (\delta_{q} (xq)^2).
\end{align*} \qed \end{prop}

This formula coincides with the expression found in \cite{KMM}, where the Frobenius-Schur indicators of  cocentral abelian extensions are computed, {\it i.e.}, extensions giving rise to the trivial action $\fde: F \times Q \to Q$. Corollary \ref{fs-chi} gives also an alternative expression in terms of the 3-cocycle $\omega$ attached to $\sigma$ and $\tau$ via the Kac exact sequence.

\subsection{Twisted quantum doubles}\label{tw-db} Let $G$ be a finite group and let $\omega$ be 3-cocycle on $G$. Consider the Dijgraaf-Pasquier-Roche quasi Hopf algebra $D^{\omega}G$, also called the \emph{twisted} quantum double of $G$ \cite{dpr}. By the results in \cite{gp-ttic}, a semisimple quasi-Hopf algebra $H$ is group theoretical if and only if its quantum double is gauge equivalent to a quasi-Hopf algebra $D^{\omega}G$. The Frobenius-Schur indicators for $D^{\omega}G$ have been computed in \cite{m-ng}, and seen to coincide in this case with the indicators introduced by Bantay \cite{bantay}.

\bigbreak
It is shown in \cite{ostrik} that the category $\Rep D^{\omega}G$ is equivalent to  $\C(G \times G,  \widetilde \omega, \Delta(G), 1)$, where $\Delta(G) \simeq G$ is the diagonal subgroup of $G \times G$, and $\widetilde \omega$ is the 3-cocycle on $G \times G$ given by $\widetilde \omega = p_1^*\omega (p_2^*\omega)^{-1}$; that is,
\begin{equation} \widetilde \omega ((a_1, a_2), (b_1,  b_2), (c_1, c_2)) = \omega(a_1, b_1, c_1)
\, \omega(a_2, b_2, c_2)^{-1},
\end{equation}
for all $a_i, b_i \in G$.

Thus our Corollary \ref{fs-chi} gives an alternative formula for the Frobenius-Schur indicators of $D^{\omega}G$ in terms of an appropriate normalization of the 3-cocycle $\widetilde \omega$.

\bigbreak
\centerline{{\sc Acknowledgement}}

\bigbreak
This work  was began during a postdoctoral stay at the Department of Mathematics of  the \' Ecole Normale Sup\' erieure, Paris. The author is grateful to Marc Rosso for his kind hospitality.
She also thanks   Peter  Schauenburg for helpful comments concerning Theorem \ref{q-antip}. Special thanks to Susan Montgomery for interesting discussions, and the Department of Mathematics of the University of Southern California for their support and  warm hospitality during her visit in October 2003.

\end{document}